\documentclass[12pt]{iopart}

\usepackage{graphicx}
\usepackage{subfigure}
\usepackage{amssymb}
\usepackage{footnote}
\usepackage{amsthm}
\usepackage{algorithmicx}
\usepackage{algorithm}
\usepackage{cases}

\newtheorem{theorem}{Theorem}
\newtheorem{lemma}[theorem]{Lemma}
\newtheorem{corollary}[theorem]{corollary}
\newtheorem{definition}{Definition}

\newcommand{\mc}[1]{{\ensuremath{\mathcal{#1}}}}

\newcommand{\mb}[1]{{\ensuremath{\mathbf{#1}}}}

\begin{document}
\title{On the inverse of the Hadamard product of a full rank matrix and an angle matrix}

\author{Yao-Jen Liang}

\address{Department of Electronic Engineering, National ILan University, Yilan 260, Taiwan (R.O.C.)}
\ead{yaojen@niu.edu.tw}
\vspace{10pt}

\begin{abstract}
By the definition of an angle matrix, we investigate the inverse of the Hadamard product of a full rank and an angle matrices.
Our proof involves standard matrix analysis.
It enriches the algebra of Hadamard products.
\end{abstract}


%
%
%

\section{Introduction}\label{sec1}
For two $m \times n$ matrices $\mb{A} = (a_{ij}), \mb{B} = (b_{ij})$, their Hadamard product (entry-wise product) is
defined to be $\mb{A}\circ \mb{B} = (a_{ij} b_{ij})$. Given an
integer $k$, the $k$th Hadamard power of $\mb{A}$ is $\mb{A}^{\circ k} = (a^k_{ij})$ \cite{Horn1991} \cite{Gregor1998}.
The Hadamard inverse of $\mb{A}$ is thus $\mb{A}^{\circ (-1)} = (a^{-1}_{ij})$, if all $a_{ij}\neq 0$.

The properties of Hadamard products and their applications have received much attention over decades as in \cite{Styan1973} \cite{Ando1979} \cite{Horn1985} \cite{Horn1991} \cite{Mathias1992} \cite{Fiedler1988} \cite{Gregor1998} \cite{Fallat2007} \cite{Liang2010} \cite{Liang2012} \cite{Yang2019} \cite{Damm2023}.
Applications of Hadamard products in artificial intelligence can be found in \cite{Zhang2020}. But to the best of our knowledge,
the inverse of
Hadamard products has not been addressed so far, although the fields of machine learning and artificial intelligence are very active now.
This motivates us to investigate the inverse of Hadamard products of some special matrices.

Throughout we consider complex matrices and denote by $\mb{A}^*$, $\mb{A}^T$, and $\mb{A}^H$ the conjugate, transpose, and
conjugate transpose of $\mb{A}$, respectively.
Further, we denote by $\mb{A}^{\circ (-T)}$ the transpose of $\mb{A}^{\circ (-1)}$ for notational simplicity.
$|\mb{A}|$ is the determinant of matrix $\mb{A}$.
The set of all complex $n\times n$ matrices is denoted as $\mc{M}_n$.

\section{Main Results}\label{sec2}
We will start with introducing the definition of a angle vector and an angle matrix with rank one.
\begin{definition}\label{def1}
Define an angle vector $\mb{v}$ and a rank-one angle matrix $\mb{\Theta}$ as\\
(a) $\mb{v} = (e^{j\theta_i})\in \mc{C}^n$ with $\theta_i \in \mc{R}, \forall i$. \\
(b) Given $\mb{v} = (e^{j\theta_i})\in \mc{C}^m$ and $\mb{u} = (e^{j\phi_i})\in \mc{C}^n$,
\begin{eqnarray}
\mb{\Theta}&=\mb{v}\cdot\mb{u}^T_n \nonumber \\
    &= \left(
         \begin{array}{ccc}
           e^{j(\theta_1+\phi_1)} & \ldots & e^{j(\theta_1+\phi_n)} \\
           \vdots & \vdots & \vdots \\
           e^{j(\theta_m+\phi_1)} & \ldots & e^{j(\theta_m+\phi_n)} \\
         \end{array}
       \right).\label{eq_Theta}
\end{eqnarray}
\end{definition}
Next, we review the following definition.
\begin{definition}(\cite{Horn1985})\label{def2}
  Let $\mb{A} = (a_{ij})\in \mc{M}_n$. Let $\mb{M}_{a_{ij}}\in \mc{M}_{n-1}$ obtained from \mb{A} by deleting the row and column containing $a_{ij}$. The cofactor $\mb{A}_{ij}$ of $a_{ij}$ is defined by $\mb{A}_{ij}=(-1)^{i+j}|\mb{M}_{a_{ji}}|$
\end{definition}

To proceed we need the following lemmas.
\begin{lemma}\label{lemma1}
  From (\ref{eq_Theta}), $\mb{\Theta}^{\circ (-T)} = \mb{\Theta}^H$.
\end{lemma}

\begin{lemma}\label{lemma2}
  The determinant of $(\mb{A}\circ \mb{\Theta})$ can be expressed as
  \begin{equation}\label{eq_det}
  |\mb{A}\circ \mb{\Theta}| = |\mb{A}|e^{j\sum_i (\theta_i+\phi_i)}.
  \end{equation}
\end{lemma}

First considering a nonsingular matrix and an angle matrix as in Definition \ref{def1}, we shall prove the following theorem. Then,
in Theorem \ref{thm2}, we will extend the result to more general full-rank matrices.
\begin{theorem}\label{thm1}
If $\mb{A}\in \mc{M}_n$ is nonsingular and $\mb{\Theta}$ is defined as in (\ref{eq_Theta}), the inverse of $(\mb{A}\circ \mb{\Theta})$ can be expressed as
\begin{equation}\label{eq_inv}
 (\mb{A}\circ \mb{\Theta})^{-1} = \mb{A}^{-1}\circ \mb{\Theta}^H.
\end{equation}
\end{theorem}

\begin{proof}
We use induction on $n$. The case $n=1$, (\ref{eq_inv}) is clearly true.

Subsequently, for $n=2$, we have
\begin{eqnarray}
\mb{A}&=\left(
                 \begin{array}{cc}
                   a_{11} & a_{12} \\
                   a_{21} & a_{22} \\
                 \end{array}
               \right),\nonumber\\
\mb{\Theta}&=\left(
              \begin{array}{cc}
                e^{j(\theta_1+\phi_1)} & e^{j(\theta_1+\phi_2)} \\
                e^{j(\theta_2+\phi_1)} & e^{j(\theta_2+\phi_2)} \\
              \end{array}
            \right),\nonumber\\
|\mb{A}\circ \mb{\Theta}| &= \left|\begin{array}{cc}
                   a_{11}e^{j(\theta_1+\phi_1)} & a_{12}e^{j(\theta_1+\phi_2)} \\
                   a_{21}e^{j(\theta_2+\phi_1)} & a_{22}e^{j(\theta_2+\phi_2)} \\
                 \end{array}\right| = |\mb{A}|e^{j(\theta_1+\theta_2+\phi_1+\phi_2)},\label{eq_det2}\\
(\mb{A}\circ \mb{\Theta})^{-1}&=\frac{1}{|\mb{A}\circ \mb{\Theta}|}\left(
                                           \begin{array}{cc}
                                             a_{22}e^{j(\theta_2+\phi_2)}  & -a_{12}e^{j(\theta_1+\phi_2)} \\
                                             -a_{21}e^{j(\theta_2+\phi_1)} & a_{11}e^{j(\theta_1+\phi_1)} \\
                                           \end{array}
                                         \right),\nonumber\\
                              &\stackrel{(\ref{eq_det2})}{=}\frac{1}{|\mb{A}|}\left(
                                           \begin{array}{cc}
                                             a_{22}  & -a_{12} \\
                                             -a_{21} & a_{11} \\
                                           \end{array}
                                         \right)\circ \left(
                                                        \begin{array}{cc}
                                                          e^{-j(\theta_1+\phi_1)} & e^{-j(\theta_2+\phi_1)} \\
                                                          e^{-j(\theta_1+\phi_2)} & e^{-j(\theta_2+\phi_2)} \\
                                                        \end{array}
                                                      \right)
                                         ,\nonumber\\
                              &=\mb{A}^{-1}\circ \mb{\Theta}^{\circ -T} = \mb{A}^{-1}\circ \mb{\Theta}^H.\nonumber
\end{eqnarray}
Then, we assume the result holds for $n=k-1$.

Finally, as $n=k$, we set $\mb{C}=\mb{B}^{-1}=(\mb{A}\circ \mb{\Theta})^{-1}$.
According to Definition \ref{def2},
it is well-known \cite{Horn1985} that
\begin{eqnarray}\label{eq_inv2}
  c_{ij} &=\frac{\mb{B}_{ij}}{|\mb{B}|} =\frac{(-1)^{i+j}|\mb{M}_{a_{ji}}\circ \mb{M}_{\theta_{ji}}|}{|\mb{A}\circ \mb{\Theta}|}
          \stackrel{(\ref{eq_det})}{=}\frac{(-1)^{i+j}|\mb{M}_{a_{ji}}|
          e^{j(\sum_{m\neq j}\theta_m+\sum_{m\neq i}\phi_m)}}{|\mb{A}|e^{j\sum_i (\theta_i+\phi_i)}}\\
         &= \frac{\mb{A}_{ij}}{|\mb{A}|}e^{-j(\theta_j+\phi_i)}.
\end{eqnarray}
Thus,
\begin{eqnarray}
\mb{C} &=(\mb{A}\circ \mb{\Theta})^{-1}=\frac{ \left(
           \begin{array}{ccc}
             \mb{A}_{11} & \ldots & \mb{A}_{1k} \\
             \vdots & \ddots & \vdots \\
             \mb{A}_{k1} & \ldots & \mb{A}_{kk} \\
           \end{array}
         \right)} {|\mb{A}|} \circ \left(
           \begin{array}{ccc}
             e^{-j(\theta_1+\phi_1)} & \ldots & e^{-j(\theta_k+\phi_1)} \\
             \vdots & \ddots & \vdots \\
             e^{-j(\theta_1+\phi_k)} & \ldots & e^{-j(\theta_k+\phi_k)} \\
           \end{array}
         \right)\nonumber \\
       &=\mb{A}^{-1}\circ \mb{\Theta}^H.
\end{eqnarray}
\end{proof}


\begin{corollary}
If $\mb{A}\in \mc{M}_n$ is nonsingular and $\mb{\Theta}$ is defined as in (\ref{eq_Theta}), $(\mb{A}\circ \mb{\Theta}^T)^{-1} = \mb{A}^{-1}\circ \mb{\Theta}^*$.
\end{corollary}
\begin{proof}
  Substituting $\mb{\Theta}$ by $\mb{\Theta}^T$, we prove this from Theorem \ref{thm1}.
\end{proof}

For the presentation of Moore-Penrose inverse of Hadamard products, we need the following lemmas.
\begin{lemma}\label{lemma3}
From Definition \ref{def1}, $\mb{\Theta}$ is as in (\ref{eq_Theta}), then \\
(a) $\mb{\Theta}^H\mb{\Theta}$ and
 $\mb{\Theta}\mb{\Theta}^H$ are $m$-times and $n$-times of some angle matrices, respectively. \\
(b) Moreover, $\mb{\Theta}^H\mb{\Theta}\mb{\Theta}^H=mn\mb{\Theta}^H$ and $\mb{\Theta}\mb{\Theta}^H\mb{\Theta}=mn\mb{\Theta}$.
\end{lemma}
\begin{proof}
Here the case of $m=n=2$ is proved, while the more general cases ($m>2, n>2$) could be proved in a similar way.
Let $\mb{v} = (e^{j\theta_1}, e^{j\theta_2})^T$ and $\mb{u} = (e^{j\phi_1}, e^{j\phi_2})^T$, we acquire
\begin{equation}\label{eq_Theta2}
  \mb{\Theta}=\left(
                \begin{array}{cc}
                  e^{j(\theta_1+\phi_1)} & e^{j(\theta_1+\phi_2)} \\
                  e^{j(\theta_2+\phi_1)} & e^{j(\theta_2+\phi_2)} \\
                \end{array}
              \right).
\end{equation}
(a)
\begin{eqnarray}\label{eq_thetatheta}
  \mb{\Theta}^H\mb{\Theta} &= \left(
                                    \begin{array}{cc}
                                      e^{-j(\theta_1+\phi_1)} & e^{-j(\theta_2+\phi_1)} \\
                                      e^{-j(\theta_1+\phi_2)} & e^{-j(\theta_2+\phi_2)} \\
                                    \end{array}
                                  \right)\left(
                                             \begin{array}{cc}
                                               e^{j(\theta_1+\phi_1)} & e^{j(\theta_1+\phi_2)} \\
                                               e^{j(\theta_2+\phi_1)} & e^{j(\theta_2+\phi_2)} \\
                                             \end{array}
                                           \right)\nonumber \\
                                &= \left(
                                    \begin{array}{cc}
                                      2 & 2e^{j(\phi_2-\phi_1)} \\
                                      2e^{j(\phi_1-\phi_2)} & 2 \\
                                    \end{array}
                                  \right) \\
&= 2\left(
                                              \begin{array}{c}
                                                1 \\
                                                e^{j(\phi_1-\phi_2)} \\
                                              \end{array}
                                            \right)\left(
                                                     \begin{array}{cc}
                                                       1 & e^{j(\phi_2-\phi_1)} \\
                                                     \end{array}
                                                   \right).\nonumber
\end{eqnarray}
Similarly,
\begin{equation*}
  \mb{\Theta}\mb{\Theta}^H = 2\left(
                                              \begin{array}{c}
                                                1 \\
                                                e^{j(\theta_2-\theta_1)} \\
                                              \end{array}
                                            \right)\left(
                                                     \begin{array}{cc}
                                                       1 & e^{j(\theta_1-\theta_2)} \\
                                                     \end{array}
                                                   \right).
\end{equation*}
(b)
\begin{eqnarray}
 \mb{\Theta}^H\mb{\Theta}\mb{\Theta}^H &\stackrel{(\ref{eq_Theta2})(\ref{eq_thetatheta})}{=}
\left(
   \begin{array}{cc}
     2 & 2e^{j(\phi_2-\phi_1)} \\
     2e^{j(\phi_1-\phi_2)} & 2 \\
   \end{array}
\right)\left(
                \begin{array}{cc}
                  e^{-j(\theta_1+\phi_1)} & e^{-j(\theta_2+\phi_1)} \\
                  e^{-j(\theta_1+\phi_2)} & e^{-j(\theta_2+\phi_2)} \\
                \end{array}
              \right)\nonumber \\
&= \left(
     \begin{array}{cc}
       4e^{-j(\theta_1+\phi_1)} & 4e^{-j(\theta_2+\phi_1)} \\
       4e^{-j(\theta_1+\phi_2)} & 4e^{-j(\theta_2+\phi_2)} \\
     \end{array}
   \right)=2\cdot 2\cdot \mb{\Theta}^H.\label{eq_lemma6}
\end{eqnarray}
Clearly, $\mb{\Theta}\mb{\Theta}^H\mb{\Theta}=( \mb{\Theta}^H\mb{\Theta}\mb{\Theta}^H)^H=2\cdot 2\cdot\mb{\Theta}$.
\end{proof}

\begin{lemma}(\cite{Penrose1955}\cite{Golub1996})\label{lemma4}
The Moore–Penrose inverse (a.k.a. pseudoinverse) $\mb{A}^+$ exists for any $m\times n$ matrix $\mb{A}$.
If furthermore $\mb{A}$ is full rank, that is, its rank is $\min(m,n)$, then\\
(a) When $\mb{A}$ has linearly independent columns (equivalently, $m > n$, $\mb{A}^H\mb{A}$ is invertible),
\begin{equation}\label{eq_linv}
  \mb{A}^+ = (\mb{A}^H\mb{A})^{-1}\mb{A}^H, \, 
\end{equation}
(b) On the other hand, if $\mb{A}$ has linearly independent rows (equivalently, $m < n$, $\mb{A}\mb{A}^H$ is invertible),
\begin{equation}\label{eq_rinv}
  \mb{A}^+ = \mb{A}^H(\mb{A}\mb{A}^H)^{-1}, \, 
\end{equation}
\end{lemma}

The pseudoinverse of the Hadamard product of a full rank matrix and an angle matrix is thus given in the following theorem.
\begin{theorem}\label{thm2}
If $m\times n$ matrix $\mb{A}$ is full rank and angle matrix $\mb{\Theta}$ is defined as in (\ref{eq_Theta}), then
\begin{subnumcases}{(\mb{A}\circ \mb{\Theta})^+=\mb{A}^+\circ \mb{\Theta}^H=}
    ((\mb{A}^H\mb{A})^{-1}\mb{A}^H) \circ \mb{\Theta}^H, & $m>n;$ \label{subeq1}\\
    \mb{A}^{-1}\circ \mb{\Theta}^H, & $m=n;$ \label{subeq2}\\
    (\mb{A}^H(\mb{A}\mb{A}^H)^{-1}) \circ \mb{\Theta}^H, & $m<n.$ \label{subeq3}
\end{subnumcases}
\end{theorem}
\begin{proof}
Let $\mb{A}=\left(
     \begin{array}{ccc}
       a_{11} & \ldots & a_{1n} \\
       \vdots & \vdots & \vdots \\
       a_{m1} & \ldots & a_{mn} \\
     \end{array}
   \right)$ and $\mb{\Theta}=\left(
                  \begin{array}{ccc}
                    e^{j(\theta_1+\phi_1)} & \ldots & e^{j(\theta_1+\phi_n)} \\
                    \vdots & \vdots & \vdots \\
                    e^{j(\theta_m+\phi_1)} & \ldots & e^{(\theta_m+\phi_n)} \\
                  \end{array}
                \right)$.
We first show that
\begin{eqnarray*}
&(\mb{A}\circ \mb{\Theta})^H(\mb{A}\circ \mb{\Theta}) \\
&=\left(
    \begin{array}{ccc}
      a^*_{11}e^{-j(\theta_1+\phi_1)} & \ldots & a^*_{m1}e^{-j(\theta_m+\phi_1)} \\
      \vdots & \vdots & \vdots \\
      a^*_{1n}e^{-j(\theta_1+\phi_n)} & \ldots & a^*_{mn}e^{-j(\theta_m+\phi_n)} \\
    \end{array}
  \right)\\
& \times \left(
    \begin{array}{ccc}
      a_{11}e^{j(\theta_1+\phi_1)} & \ldots & a_{1n}e^{j(\theta_1+\phi_n)} \\
      \vdots & \vdots & \vdots \\
      a_{m1}e^{j(\theta_m+\phi_1)} & \ldots & a_{mn}e^{j(\theta_m+\phi_n)} \\
    \end{array}
  \right)\\
&=\left(
    \begin{array}{ccc}
      \sum_l|a_{l1}|^2 & \ldots & \sum_la^*_{l1}a_{ln}e^{j(\phi_n-\phi_1)} \\
      \vdots & \ddots & \vdots \\
      \sum_la^*_{ln}a_{l1}e^{j(\phi_1-\phi_n)} & \ldots & \sum_l|a_{ln}|^2 \\
    \end{array}
  \right)\\
&=\left(
    \begin{array}{ccc}
      \sum_l|a_{l1}|^2 & \ldots & \sum_la^*_{l1}a_{ln} \\
      \vdots & \ddots & \vdots \\
      \sum_la^*_{ln}a_{l1} & \ldots & \sum_l|a_{ln}|^2 \\
    \end{array}
  \right)\circ \left(
    \begin{array}{ccc}
      1 & \ldots & e^{j(\theta_n-\theta_1)} \\
      \vdots & \ddots & \vdots \\
      e^{j(\phi_1-\phi_n)} & \ldots & 1 \\
    \end{array}
  \right)\\
&\stackrel{(\ref{eq_thetatheta})}=(\mb{A}^H\mb{A})\circ(\frac{1}{m}\mb{\Theta}^H\mb{\Theta}). 
\end{eqnarray*}
Subsequently, it suffices to prove the case as in (\ref{subeq1}), since (\ref{subeq2}) is from Theorem \ref{thm1} and (\ref{subeq3}) can be derived similar to the derivation of (\ref{subeq1}).
When $m>n$,
\begin{eqnarray}
(\mb{A}\circ \mb{\Theta})^+ &\stackrel{(\ref{eq_linv})}=
[(\mb{A}\circ \mb{\Theta})^H(\mb{A}\circ \mb{\Theta})]^{-1}(\mb{A}\circ \mb{\Theta})^H \\
&= [(\mb{A}^H\mb{A})\circ (\frac{1}{m}\mb{\Theta}^H\mb{\Theta})]^{-1} (\mb{A}^H\circ \mb{\Theta}^H)\\
&\stackrel{(\ref{eq_inv})}{=}m[(\mb{A}^H\mb{A})^{-1}\circ (\mb{\Theta}^H\mb{\Theta})^H] (\mb{A}^H\circ \mb{\Theta}^H)\\
&= m((\mb{A}^H\mb{A})^{-1}\mb{A}^H) \circ(\frac{1}{m}\mb{\Theta}^H\mb{\Theta}\mb{\Theta}^H)\\
&\stackrel{(\ref{eq_lemma6})}{=} ((\mb{A}^H\mb{A})^{-1}\mb{A}^H) \circ \mb{\Theta}^H\\
&=\mb{A}^+\circ \mb{\Theta}^H.
\end{eqnarray}
\end{proof}

\end{document}